# A Conjectural Description of the Tautological Ring of the Moduli Space of Curves

Carel Faber

The purpose of this paper is to formulate a number of conjectures giving a rather complete description of the tautological ring of $\mathcal{M}_g$ and to discuss the evidence for these conjectures.

Denote by $\mathcal{M}_g$ the moduli space of smooth curves of genus $g \geq 2$ over an algebraically closed field (of arbitrary characteristic). For every $m \geq 3$ that is not divisible by the characteristic, this space is the quotient by a finite group of the nonsingular moduli space of smooth curves with a symplectic level-$m$ structure. Hence the Chow ring $A^*(\mathcal{M}_g)$ can be defined easily, since we will be working with $\mathbb{Q}$-coefficients throughout this paper. (Cf. [Fu], Example 8.3.12.)

Some other relevant spaces are the moduli spaces $\mathcal{M}_{g,n}$ of smooth $n$-pointed curves $(C; x_1, \ldots, x_n)$ of genus $g$, with $x_i \neq x_j$ for $i \neq j$, defined whenever $2g-2+n > 0$, and, denoting $\mathcal{M}_{g,1}$ by $\mathcal{C}_g$, the spaces $\mathcal{C}_g^n$, the $n$-fold fibre products of $\mathcal{C}_g$ over $\mathcal{M}_g$, parametrizing smooth curves of genus $g$ with $n$-tuples of not necessarily distinct points. For all these spaces, the Chow ring can be defined as above.

Between these spaces, there are many natural morphisms forgetting one or more points; these will usually be denoted $\pi$. The most important of these is $\pi : \mathcal{C}_g \to \mathcal{M}_g$; its relative dualizing sheaf $\omega_\pi$ will often be denoted by $\omega$. This is a $Q$-line bundle on $\mathcal{C}_g$, cf. [Mu], p. 299. Writing $K = c_1(\omega) \in A^1(\mathcal{C}_g)$ we define (following Mumford)

$$\kappa_i := \pi_*(K^{i+1}) \in A^i(\mathcal{M}_g)$$

using the ring structure in $A^*(\mathcal{C}_g)$ and the proper push-forward for Chow groups. Note that $\kappa_0 = 2g-2$ and $\kappa_{-1} = 0$.

Another way to produce natural classes in $A^*(\mathcal{M}_g)$ is to use the Hodge bundle $\mathbb{E} = \pi_*\omega$, a locally free $Q$-sheaf of rank $g$ on $\mathcal{M}_g$. It is the pull-back of a bundle on $\mathcal{A}_g$, the moduli space of principally polarized abelian varieties of dimension $g$, via the morphism $t : \mathcal{M}_g \to \mathcal{A}_g$ sending a curve to its Jacobian. We follow Mumford in writing

$$\lambda_i := c_i(\mathbb{E}) \in A^i(\mathcal{M}_g).$$

(Thus $\lambda_0 = 1$ and $\lambda_i = 0$ for $i > g$; the divisor class $\lambda_1$ is often denoted $\lambda$.)

The $\kappa_i$ and $\lambda_i$ are called the *tautological classes*; we define the *tautological ring* of $\mathcal{M}_g$ to be the $\mathbb{Q}$-subalgebra of $A^*(\mathcal{M}_g)$ generated by the tautological classes $\kappa_i$ and $\lambda_i$. We denote this ring by $R^*(\mathcal{M}_g)$.

It can be shown that the classes of many geometrically defined subvarieties of $\mathcal{M}_g$ lie in $R^*(\mathcal{M}_g)$. Examples of subvarieties defined in terms of linear systems for which this happens, will play an important role in this paper. Mumford gives examples (in §7 of [Mu]) of subvarieties parametrizing curves with special Weierstrass points; these are actually special cases of subvarieties defined in terms of linear systems with certain prescribed types of ramification. Working over $\mathbb{C}$, we have Harer's result ([Ha 1]) that $H^2(\mathcal{M}_g)$ is one-dimensional (for $g \geq 3$), implying that

$$A^1(\mathcal{M}_g) = R^1(\mathcal{M}_g) \cong \mathbb{Q}$$

for $g \geq 3$ (note that $\kappa_1 = 12\lambda$, cf. [Mu] p. 306). (Unfortunately, an algebraic proof of Harer's result is not known.)

For $g \leq 5$, the tautological ring $R^*(\mathcal{M}_g)$ is equal to the Chow ring $A^*(\mathcal{M}_g)$, in characteristic 0. This was shown by Mumford for $g = 2$ (cf. [Mu] p. 318; it is an immediate consequence of Igusa's description of $\mathcal{M}_2$), by the author for genera 3 and 4 ([Fa 1], [Fa 2]) and by Izadi for $g = 5$ ([Iz]). However, it doesn't seem possible that this hold for all $g$. The idea is that the recent result of Pikaart that the cohomology of $\overline{\mathcal{M}}_g$ is not of Tate type for $g$ large ([Pi], Cor. 4.7) should imply that the Chow groups of $\overline{\mathcal{M}}_g$ over $\mathbb{C}$ don't map injectively to the cohomology groups; here I am assuming that the result of Jannsen for smooth projective varieties over a universal domain ([Ja], Thm. 3.6(a)) can be extended to varieties like $\overline{\mathcal{M}}_g$ over $\mathbb{C}$. It would seem then that a similar result would hold for the Chow ring of $\mathcal{M}_g$. The question whether $R^*(\mathcal{M}_g)$ and $A^*(\mathcal{M}_g)$ have the same image in $H^*(\mathcal{M}_g)$ appears to be open.

We now formulate some known results about the tautological ring of $\mathcal{M}_g$. First of all, Mumford shows ([Mu], §§5, 6) that the ring $R^*(\mathcal{M}_g)$ is generated by the $g-2$ classes $\kappa_1, \ldots, \kappa_{g-2}$. The proof has 2 ingredients:

a. Applying the Grothendieck-Riemann-Roch theorem to $\pi : \mathcal{C}_g \to \mathcal{M}_g$ and $\omega_\pi$ gives an expression for the Chern character of the Hodge bundle $\mathbb{E}$ in terms of the $\kappa_i$. The resulting expressions for the Chern classes $\lambda_i$ of $\mathbb{E}$ in the $\kappa_i$ can be concisely formulated as an identity of formal power series in $t$:

$$\sum_{i=0}^{\infty} \lambda_i t^i = \exp\left(\sum_{i=1}^{\infty} \frac{B_{2i}\,\kappa_{2i-1}}{2i(2i-1)} t^{2i-1}\right).$$

Here $B_{2i}$ are the Bernoulli numbers with signs ($B_2 = 1/6$, $B_4 = -1/30$, ...). For example

$$\lambda_1 = \frac{1}{12}\kappa_1, \qquad \lambda_2 = \frac{1}{2}\lambda_1^2 = \frac{1}{288}\kappa_1^2, \qquad \lambda_3 = \frac{1}{6}\left(\frac{\kappa_1}{12}\right)^3 - \frac{1}{360}\kappa_3.$$

So all the $\lambda_i$ can be expressed in the odd $\kappa_i$. (Also, the odd $\kappa_i$ with $i > g$ can be expressed in the lower (odd) kappa's; this will not be used in the proof that $\kappa_1, \ldots, \kappa_{g-2}$ generate, therefore it gives relations between the latter classes in odd degrees greater than $g$.)

b. On a non-singular curve, the relative dualizing sheaf is generated by its global sections. This may be formulated universally as the surjectivity of the natural map $\pi^*\mathbb{E} \to \omega$ of locally free sheaves on $\mathcal{C}_g$. The kernel is then locally free of rank $g-1$ so that its Chern classes vanish in degrees greater than $g-1$. Hence

$$c_j(\pi^*\mathbb{E} - \omega) = 0 \qquad \forall j \geq g\,,$$

the difference being taken in the Grothendieck group. Pushing-down to $\mathcal{M}_g$ gives relations between the lambda's and the kappa's in every degree $\geq g-1$.

To obtain the desired result, one needs to check that in degrees $g-1$ and $g$ the 2 relations are independent. For this, Mumford uses an estimate on the size of the Bernoulli numbers; alternatively, one can use (easy) congruence properties of these numbers.

From now on, when talking about relations in the tautological ring, we will mean relations between the kappa's.

Recently, Looijenga proved a strong vanishing result about the tautological ring ([Lo]):

**Theorem 1 (Looijenga).** *$R^j(\mathcal{M}_g) = 0$ for all $j > g-2$ and $R^{g-2}(\mathcal{M}_g)$ is at most one-dimensional, generated by the class of the hyperelliptic locus.*

In fact he proved a similar statement for the tautological ring of $\mathcal{C}_g^n$. This is the subring of $A^*(\mathcal{C}_g^n)$ generated by the divisor classes $K_i := pr_i^* K$ and $D_{ij}$ (the class of the diagonal $x_i = x_j$) and the pull-backs from $\mathcal{M}_g$ of the $\kappa_i$. The result is that it vanishes in degrees greater than $g-2+n$ and that it is at most one-dimensional in degree $g-2+n$, generated by the class of the locus

$$\mathcal{H}_g^n = \{(C; x_1, \ldots, x_n) : C \text{ hyperelliptic; } x_1 = \ldots = x_n = x, \text{ a Weierstrass point}\}.$$

Looijenga also gives a description of the degree $d$ part of the tautological ring of $\mathcal{C}_g^n$, for all $d$ and $n$.

This establishes an important part of one of the conjectures to be discussed in this paper. That conjecture immediately implies Diaz's theorem ([Di 2]), which gives the upper bound $g-2$ for the dimension of a complete subvariety of $\mathcal{M}_g$ (in char. 0). Diaz used a flag of subvarieties of $\mathcal{M}_g$ that refines the flag introduced by Arbarello ([Ar]). By using another refinement of Arbarello's flag and by simplifying parts of Diaz's work, Looijenga is able to express every tautological class of degree $d$ as a linear combination of the classes of the irreducible components of a specific geometrically defined locus of $n$-pointed curves. For $d > g-2+n$, the locus is empty; hence the vanishing. In degree $g-2+n$, the locus is not irreducible. After having described the irreducible components, Looijenga proves that their classes are proportional to that of $\mathcal{H}_g^n$ by invoking the Fourier transform for abelian varieties (work of Mukai, Beauville and Deninger-Murre). As a corollary of Looijenga's theorem one obtains Diaz's result in arbitrary characteristic.

The question whether the class of the hyperelliptic locus $\mathcal{H}_g$ is actually non-zero was left open; this had been established only for $g = 3$ (due to the existence of complete curves in $\mathcal{M}_3$) and $g = 4$ (by means of a long calculation with 'test surfaces' in $\overline{\mathcal{M}}_4$ (see [Fa 2])). Note that the vanishing of $[\mathcal{H}_g]$ would imply an improvement of Diaz's bound; or conversely, the existence of a complete subvariety of dimension $g-2$ of $\mathcal{M}_g$ would imply the non-vanishing of $\kappa_1^{g-2}$ (since $\kappa_1$ is ample), hence that of $[\mathcal{H}_g]$. However, it is not even known whether $\mathcal{M}_4$ contains a complete surface.

Happily enough, we don't need the existence result for complete subvarieties to settle the non-vanishing in the top degree:

**Theorem 2.** *$\kappa_{g-2} \neq 0$ on $\mathcal{M}_g$. Hence $R^{g-2}(\mathcal{M}_g)$ is one-dimensional.*

So the classes $[\mathcal{H}_g]$ and $[\mathcal{H}_g^n]$ are non-zero as well. Before discussing the proof, we remind the reader of the fact that the classes $\kappa_i$ and $\lambda_i$, which so far have been defined only on $\mathcal{M}_g$, can be defined exactly as above on the Deligne-Mumford compactification $\overline{\mathcal{M}}_g$. See [Mu], §4.

The proof consists of two parts:

a. The class $\lambda_{g-1}\lambda_g$ vanishes on the boundary $\overline{\mathcal{M}}_g - \mathcal{M}_g$ of the moduli space.

This is a simple observation. In fact, over $\Delta_0$, the closure of the locus of irreducible singular curves, the class $\lambda_g$ already vanishes: when pulled back to $\overline{\mathcal{M}}_{g-1,2}$, the Hodge bundle on $\Delta_0$ becomes an extension of a trivial line bundle by the Hodge bundle in genus $g-1$, so its top Chern class vanishes. Over a boundary component $\Delta_i$, with $1 \leq i \leq [g/2]$, the closure of the locus of reducible singular curves consisting of one component of genus $i$ and one of genus $g-i$, the Hodge bundle becomes the direct sum of the Hodge bundles in genera $i$ and $g-i$. Now use the identity $\lambda_h^2 = 0$, valid in arbitrary genus $h$, to conclude the vanishing of $\lambda_{g-1}\lambda_g$ over $\Delta_i$.

An equivalent formulation is:

a$'$. The class $ch_{2g-1}(\mathbb{E})$ vanishes on the boundary $\overline{\mathcal{M}}_g - \mathcal{M}_g$ of the moduli space.

One may prove this directly, using on the one hand the additivity of the Chern character in exact sequences and on the other hand the vanishing of *all* the components of degree $\geq 2h$ of the Chern character of the Hodge bundle in genus $h$, an easy consequence of the vanishing of the even components proved by Mumford in [Mu], §5. Or one uses the identity

$$\lambda_{g-1}\lambda_g = (-1)^{g-1}(2g-1)! \cdot ch_{2g-1}(\mathbb{E}),$$

another consequence of Mumford's result.

b. On $\overline{\mathcal{M}}_g$ the following identity holds:

$$\kappa_{g-2}\lambda_{g-1}\lambda_g = \frac{|B_{2g}|(g-1)!}{2^g(2g)!} \quad . \tag{1}$$

This is an identity of intersection numbers (more precisely, the number on the right is the degree of the zero cycle on the left). As $B_{2g}$ doesn't vanish, this proves the theorem.

The first step in the proof of (1) is to use (at last!) the full force of Mumford's result in §5 of [Mu]: his expression for the Chern character of the Hodge bundle on $\overline{\mathcal{M}}_g$, derived by using the Grothendieck-Riemann-Roch theorem twice. Using this, (1) translates into the following identity of intersection numbers of Witten's tau-classes ([Wi]):

$$\frac{g!}{2^{g-1}(2g)!} = \langle\tau_{g-1}\tau_{2g}\rangle - \langle\tau_{3g-2}\rangle + \frac{1}{2}\sum_{j=0}^{2g-2}(-1)^j\langle\tau_{2g-2-j}\tau_j\tau_{g-1}\rangle$$
$$+\frac{1}{2}\sum_{h=1}^{g-1}\left((-1)^{g-h}\langle\tau_{3h-g}\tau_{g-1}\rangle\langle\tau_{3(g-h)-2}\rangle + (-1)^h\langle\tau_{3h-2}\rangle\langle\tau_{3(g-h)-g}\tau_{g-1}\rangle\right).$$

The second step is to invoke the Witten conjecture, proven by Kontsevich ([Wi], [Ko]). This gives a recipe to compute all intersection numbers of tau-classes: a generating function encoding all these numbers satisfies the Korteweg-de Vries equations. Together with the so-called string and dilaton equations, this determines these numbers recursively.

Proving an explicit identity as the one above, is however not necessarily straightforward from the said recipe. Define the *n-point function* as the following formal power series in $n$ variables $x_1, \ldots, x_n$ :

$$\langle \tau(x_1) \cdots \tau(x_n) \rangle = \sum_{a_1, \ldots, a_n \geq 0} \langle \tau_{a_1} \cdots \tau_{a_n} \rangle x_1^{a_1} \cdots x_n^{a_n} ,$$

which encodes all intersection numbers of $n$ tau-classes. To prove the identity above, it suffices to know 2 special 3-point functions explicitly, to wit $\langle \tau_0 \tau(x) \tau(y) \rangle$ and $\langle \tau(x) \tau(y) \tau(-y) \rangle$. The KdV-equations in the form given by Witten ([Wi], (2.33)) can be translated into simple differential equations for these functions; an initial condition is provided by the identity

$$\langle \tau_0 \tau_0 \tau(x) \rangle = \exp\left( \frac{x^3}{24} \right) ,$$

which is an immediate consequence of the Witten conjecture. In this way, the 2 special 3-point functions can be determined easily; for instance,

$$\langle \tau_0 \tau(w) \tau(z) \rangle = \exp\left( \frac{(w^3 + z^3)}{24} \right) \sum_{n \geq 0} \frac{n!}{(2n+1)!} \left( \tfrac{1}{2} wz(w+z) \right)^n ,$$

a formula we learned from Dijkgraaf ([Dij]). This finishes the (sketch of the) proof of the theorem. (Very recently, Zagier determined the general 3-point function explicitly.)

We now formulate the first conjecture about the tautological ring $R^*(\mathcal{M}_g)$. We choose to state it first in the form in which it was discussed at several occasions, as early as Spring 1993, in particular, before the 2 theorems above were proved.

**Conjecture 1.**

a. *The tautological ring $R^*(\mathcal{M}_g)$ is Gorenstein with socle in degree $g-2$. I.e., it vanishes in degrees $> g-2$, is 1-dimensional in degree $g-2$ and, when an isomorphism $R^{g-2}(\mathcal{M}_g) = \mathbb{Q}$ is fixed, the natural pairing*

$$R^i(\mathcal{M}_g) \times R^{g-2-i}(\mathcal{M}_g) \to R^{g-2}(\mathcal{M}_g) = \mathbb{Q}$$

*is perfect.*

b. *The $[g/3]$ classes $\kappa_1, \ldots, \kappa_{[g/3]}$ generate the ring, with no relations in degrees $\leq [g/3]$.*

c. *There exist explicit formulas for the proportionalities in degree $g-2$, which may be given as follows. We define expressions $\langle \tau_{d_1+1} \tau_{d_2+1} \cdots \tau_{d_k+1} \rangle$, elements of $R^{g-2}(\mathcal{M}_g)$, in 2 ways, for every partition of $g-2$ into positive integers $d_1, d_2, \ldots, d_k$ ; this allows to express every monomial $\kappa_I$ of degree $g-2$ (where $I$ is a multi-index) as a multiple of $\kappa_{g-2}$ .*

*(1)*

$$\langle \tau_{d_1+1}\tau_{d_2+1}\cdots\tau_{d_k+1}\rangle = \frac{(2g-3+k)!(2g-1)!!}{(2g-1)!\prod_{j=1}^{k}(2d_j+1)!!}\,\kappa_{g-2}\,.$$

*(2)*

$$\langle \tau_{d_1+1}\tau_{d_2+1}\cdots\tau_{d_k+1}\rangle = \sum_{\sigma\in\mathfrak{S}_k}\kappa_\sigma\,,$$

*where* $\kappa_\sigma = \kappa_{|\alpha_1|}\kappa_{|\alpha_2|}\cdots\kappa_{|\alpha_{\nu(\sigma)}|}$ *for a decomposition* $\sigma = \alpha_1\alpha_2\cdots\alpha_{\nu(\sigma)}$ *of the permutation* $\sigma$ *in disjoint cycles, including the 1-cycles; finally,* $|\alpha|$ *is defined as the sum of the elements in the cycle* $\alpha$*, where we think of* $\mathfrak{S}_k$ *as acting on the* $k$*-tuples with entries* $d_1, d_2, \ldots, d_k$ *.*

Here $(2a-1)!!$ is shorthand for $(2a)!/(2^a a!)$. A few examples will clarify the recipe given in (c) above:

$$\langle\tau_{g-1}\rangle = \kappa_{g-2}\,;$$

$$\langle\tau_i\tau_{g-i}\rangle = \kappa_{i-1}\kappa_{g-i-1} + \kappa_{g-2} = \frac{(2g-1)!!}{(2i-1)!!(2g-2i-1)!!}\,\kappa_{g-2}\,;$$

$$\langle\tau_{i+1}\tau_{j+1}\tau_{k+1}\rangle = \kappa_i\kappa_j\kappa_k + \kappa_{i+j}\kappa_k + \kappa_{i+k}\kappa_j + \kappa_{j+k}\kappa_i + 2\kappa_{g-2} \qquad (i+j+k=g-2).$$

Let me point out here that the inspiration to look at the *sums* of 'intersection numbers' (multiples of $\kappa_{g-2}$) occurring in (2) above, instead of at the numbers themselves, came entirely from the Witten conjecture ([Wi], see also [Ho]), as the notation suggests. A direct link with the actual intersection numbers of Witten's tau-classes on the compactified moduli spaces $\overline{\mathcal{M}}_{g,n}$ was not available at the time, however; it is now, via the class $ch_{2g-1}(\mathbb{E})$ mentioned in the sketch of the proof of Theorem 2. The resulting conjectural identity between the latter numbers is:

$$\frac{(2g-3+k)!}{2^{2g-1}(2g-1)!}\cdot\frac{1}{\prod_{j=1}^{k}(2e_j-1)!!} = \langle\tau_{e_1}\cdots\tau_{e_k}\tau_{2g}\rangle - \sum_{j=1}^{k}\langle\tau_{e_1}\cdots\tau_{e_{j-1}}\tau_{e_j+2g-1}\tau_{e_{j+1}}\cdots\tau_{e_k}\rangle$$
$$+\frac{1}{2}\sum_{j=0}^{2g-2}(-1)^j\langle\tau_{2g-2-j}\tau_j\tau_{e_1}\cdots\tau_{e_k}\rangle + \frac{1}{2}\sum_{\underline{k}=I\amalg J}\sum_{j=0}^{2g-2}(-1)^j\langle\tau_j\prod_{i\in I}\tau_{e_i}\rangle\langle\tau_{2g-2-j}\prod_{i\in J}\tau_{e_i}\rangle$$

where $\underline{k} = \{1,2,\ldots,k\}$ and $\sum_{j=1}^{k}(e_j-1) = g-2$. (We proved that it is compatible with the string and dilaton equations, so $e_j \geq 2$ may be assumed.)

Returning to $\mathcal{M}_g$, we can give the proportionality factor for $\kappa_1^{g-2}$ explicitly:

$$\kappa_1^{g-2} = \frac{1}{g-1}\,2^{2g-5}\big((g-2)!\big)^2\kappa_{g-2}$$

a formula which we had observed 'experimentally' and which was proven instantaneously by Zagier from (c) above.

Note that parts (a) and (c) of the conjecture implicitly determine the dimension of the $\mathbb{Q}$-vector space $R^i(\mathcal{M}_g)$: it is the rank of the $p(i)$ by $p(g-2-i)$ matrix (with $p$ the partition function) whose entries are the 'intersection numbers' $r_{IJ}$ of monomials $\kappa_I$ of degree $i$ and $\kappa_J$ of degree $g-2-i$ given by $\kappa_I\kappa_J = r_{IJ}\kappa_{g-2}$. So the second half of part (b) of the conjecture is the claim that this matrix is of maximal rank whenever $3i \leq g$. Unfortunately, we have not been able to derive an explicit formula for the dimension of $R^i(\mathcal{M}_g)$ in this manner. In joint work with Zagier, we found a relatively simple formula that fits with the data obtained sofar from computations; this will be discussed later.

Part (a) of the conjecture may be rephrased to say that $R^*(\mathcal{M}_g)$ has the Poincaré Duality property enjoyed by the ring of algebraic cohomology classes (with $\mathbb{Q}$-coefficients) of a nonsingular projective variety of dimension $g-2$. In light of this, Thaddeus asked the question whether $R^*(\mathcal{M}_g)$ also satisfies the other properties such a ring is known (resp. conjectured) to have in char. 0 (resp. in char. $p > 0$). (Cf. Grothendieck's paper [Gr] for a discussion of these.) After having examined the available evidence, we feel confident enough to extend Conjecture 1:

**Conjecture 1(bis).** *In addition to the properties mentioned in Conjecture 1, $R^*(\mathcal{M}_g)$ 'behaves like' the algebraic cohomology ring of a nonsingular projective variety of dimension $g-2$; i.e., it satisfies the Hard Lefschetz and Hodge Positivity properties with respect to the class $\kappa_1$.*

We don't have a particular candidate for such a projective variety. Diaz's upper bound allows for the existence of such a variety lying inside $\mathcal{M}_g$, although presumably it will have at least quotient singularities in that case. For a brief discussion of the relation between the conjectured form of the tautological ring and the occurrence of complete subvarieties inside moduli space, see the concluding remarks.

So as not to lose the interest of the skeptical reader, we state the following result.

**Theorem 3.** *Conjectures 1 and 1(bis) are true for all $g \leq 15$.*

In order to explain how we could settle these conjectures for the values of $g$ mentioned, we introduce certain sheaves on the spaces $\mathcal{C}_g^d$, the $d$-fold fibre products of the universal curve $\mathcal{C}_g$ over $\mathcal{M}_g$.

Consider the projection $\pi = \pi_{\{1,\ldots,d\}} : \mathcal{C}_g^{d+1} \to \mathcal{C}_g^d$ that forgets the $(d+1)$-st point. Denote by $\Delta_{d+1}$ the sum of the $d$ divisors $D_{1,d+1}, \ldots, D_{d,d+1}$ as well as (by abuse of notation) its class:

$$\Delta_{d+1} = D_{1,d+1} + \ldots + D_{d,d+1}\,.$$

Further, denote by $\omega_i$ the $(Q$-)line bundle on $\mathcal{C}_g^n$ obtained by pulling back $\omega$ on $\mathcal{C}_g$ along the projection onto the $i$-th factor and denote its class in the codimension-1 Chow group by $K_i$.

We define a coherent sheaf $\mathbb{F}_d$ on $\mathcal{C}_g^d$ by the formula

$$\mathbb{F}_d = \pi_*(\mathcal{O}_{\Delta_{d+1}} \otimes \omega_{d+1}).$$

The sheaf $\mathbb{F}_d$ is locally free of rank $d$; its fibre at a point $(C; x_1, \ldots, x_d) = (C; D)$ is the vector space

$$H^0(C, K/K(-D)).$$

We think of $\mathbb{F}_d$ as a universal $d$-pointed jet bundle. It is invariant for the action of $\mathfrak{S}_d$. Its total Chern class can be expressed in terms of the tautological divisor classes $K_i$ and $D_{ij}$ on $\mathcal{C}_g^d$:

$$\begin{aligned} c(\mathbb{F}_d) &= (1+K_1)(1+K_2-\Delta_2)(1+K_3-\Delta_3)\cdots(1+K_d-\Delta_d) \\ &= (1+K_1)(1+K_2-D_{12})(1+K_3-D_{13}-D_{23})\cdots(1+K_d-D_{1d}\ldots-D_{d-1,d}). \end{aligned}$$

This can be proved for instance using the Grothendieck-Riemann-Roch theorem and a natural filtration on the sheaves $\mathbb{F}_n$ . We refer to [Fa 4] for the details.

The natural evaluation map of locally free sheaves on $\mathcal{C}_g^d$:

$$\varphi_d : \mathbb{E} \to \mathbb{F}_d$$

between the pulled-back Hodge bundle of rank $g$ and the bundle $\mathbb{F}_d$ of rank $d$ will be our main tool in constructing relations between tautological classes. Fibrewise the kernel over $(C;D)$ is the vector space $H^0(C,K(-D))$, whose dimension may vary with $D$.

Observe that the locus $\{\operatorname{rank}\varphi_d \le d-r\}$ parametrizes the pairs $(C;D)$ for which dim $H^0(C,K(-D)) \ge g-d+r$, equivalently, dim $H^0(C,D) \ge r+1$, in other words, for which the complete linear system $|D|$ has dimension at least $r$. So the image of this locus in $\mathcal{M}_g$ parametrizes the curves possessing a $g^r_d$. The expected codimension of the locus $\{\operatorname{rank}\varphi_d \le d-r\}$ in $\mathcal{C}_g^d$ is $r(g-d+r)$; the expected fibre dimension of the map to $\mathcal{M}_g$ is $r$; so the expected codimension in $\mathcal{M}_g$ of the locus of curves possessing a $g^r_d$ is

$$\rho = r(g-d+r) - d + r = (r+1)(g-d+r) - g,$$

the Brill-Noether number.

Porteous's formula (cf. [ACGH], [Fu]) computes the class of the locus $\{\operatorname{rank}\varphi_d \le d-r\}$ if it is either empty or has the expected codimension. The formula is:

$$\text{class } (\{\operatorname{rank}\varphi_d \le d-r\}) = \Delta_{r,g-d+r}(c(\mathbb{F}_d - \mathbb{E})).$$

Here the difference is taken in the Grothendieck group; $c(\mathbb{F}_d - \mathbb{E})$ is the formal power series in $t$ obtained as the quotient $c(\mathbb{F}_d)/c(\mathbb{E})$ of total Chern classes, this time written as polynomials in $t$. Finally,

$$\Delta_{p,q}\left(\sum_{i=0}^{\infty} c_i t^i\right) = \begin{vmatrix} c_p & c_{p+1} & \cdots & c_{p+q-1} \\ c_{p-1} & c_p & \cdots & c_{p+q-2} \\ \vdots & \vdots & \ddots & \vdots \\ c_{p-q+1} & c_{p-q+2} & \cdots & c_p \end{vmatrix}.$$

As a very important example of the above, consider the curves $C$ of genus $g$ with a $g^g_{2g-1}$. Divisors of degree $-1$ on a curve are not effective; dually this says that no curve has a $g^g_{2g-1}$. Hence

$$\{\operatorname{rank}\varphi_{2g-1} \le g-1\}$$

is a fancy way to denote the empty set. Porteous's formula applies and we find:

**Proposition 1.**

$$c_g(\mathbb{F}_{2g-1} - \mathbb{E}) = 0.$$

As the Chern classes of both $\mathbb{E}$ and $\mathbb{F}_{2g-1}$ are expressed in terms of the tautological classes, this gives a relation between tautological classes on $\mathcal{C}_g^{2g-1}$. In principle, pushing-down this relation all the way to $\mathcal{M}_g$ will give a relation between the tautological classes $\kappa_i$ and $\lambda_i$, hence between the $\kappa_i$ themselves. (One has to note that a monomial in the classes $K_i$ and $D_{ij}$ pushes down to a monomial in the $\kappa_i$; the formulas governing this will be discussed shortly.)

Unfortunately, for trivial reasons the obtained relation is identically zero: the relation of the proposition lives in codimension $g$ on $\mathcal{C}_g^{2g-1}$, hence ends up in negative codimension once pushed down to $\mathcal{M}_g$. (This might be one of the reasons why this relation apparently was not considered before.)

Fortunately however, "once zero, always zero": multiplying the relation with an arbitrary class gives another relation; in particular, multiplying it with a monomial in tautological classes gives another relation between tautological classes; after pushing-down to $\mathcal{M}_g$ we obtain relations between the $\kappa_i$ which do not obviously vanish.

Before discussing these relations, we state a variant of the proposition above:

**Proposition 1(bis).**

$$\textit{For all } d \geq 2g-1, \textit{ for all } j \geq d-g+1, \qquad c_j(\mathbb{F}_d - \mathbb{E}) = 0.$$

This is because the locus $\{\operatorname{rank} \varphi_d \leq g-1\}$ is empty; so $\varphi_d$ is an injective map of vector bundles, whence the cokernel is locally free of rank $d-g$, whence the result. (To ease the exposition, we will only use $d = 2g-1$ in the sequel.)

The available evidence suggests that the relations between the $\kappa_i$ obtained from the relations just stated by multiplying with a monomial in the $K_i$ and $D_{ij}$ and pushing-down to $\mathcal{M}_g$ are very non-trivial indeed. In fact, calculations we have done show that for $g \leq 15$ these relations generate the entire ideal of relations in the tautological ring. I.e., dividing out the polynomial ring $\mathbb{Q}[\kappa_1, \ldots, \kappa_{g-2}]$ by the ideal of relations so obtained gives a quotient ring that surjects onto the tautological ring; the quotient ring is Gorenstein with socle in degree $g-2$; because $R^{g-2}(\mathcal{M}_g)$ is non-zero by Theorem 2, the surjection is in fact an isomorphism. In this way one proves Theorem 3. Below we discuss the calculations in some detail. Because we see no reason whatsoever why the result of the calculations would be different for higher genera, we put forward the following conjecture:

**Conjecture 2.** *Let $I_g$ be the ideal of relations in the polynomial ring $\mathbb{Q}[\kappa_1, \ldots, \kappa_{g-2}]$ generated by the relations of the form*

$$\pi_*\big(M \cdot c_j(\mathbb{F}_{2g-1} - \mathbb{E})\big),$$

*with $j \geq g$ and $M$ a monomial in the $K_i$ and $D_{ij}$ and $\pi : \mathcal{C}_g^{2g-1} \to \mathcal{M}_g$ the map forgetting all the points. Then the quotient ring $\mathbb{Q}[\kappa_1, \ldots, \kappa_{g-2}]/I_g$ is Gorenstein with socle in degree $g-2$; hence it is isomorphic to the tautological ring $R^*(\mathcal{M}_g)$.*

The implied isomorphism follows from Theorem 2. As mentioned, Conjecture 2 is proved for all $g \leq 15$. (As it turns out, monomials in the $D_{ij}$ suffice.)

We now discuss the 'mechanics' of the calculation: how a relation of the form stated in Conjecture 2 actually produces a relation between the $\kappa_i$ on $\mathcal{M}_g$. Observe that the expression

$$M \cdot c_j(\mathbb{F}_{2g-1} - \mathbb{E})$$

may be expanded into a polynomial in the classes $K_i$, $D_{ij}$ and $\lambda_i$. The latter classes are pull-backs from $\mathcal{M}_g$; we will essentially always suppress this in the notation. Also, although strictly speaking we don't need it at this point, we point out that it is easy to invert $c(\mathbb{E})$:

$$c(\mathbb{E})^{-1} = c(\mathbb{E}^\vee) = 1 - \lambda_1 + \lambda_2 - \lambda_3 + \ldots + (-1)^g \lambda_g\,.$$

Hence it suffices to explain how to compute $\pi_*$ of a monomial in the classes just mentioned. The projection formula tells us that it suffices to do this for monomials in the classes $K_i$ and $D_{ij}$. Now the map $\pi$ is the composition of morphisms $\pi_d$ from $\mathcal{C}_g^d$ to $\mathcal{C}_g^{d-1}$ forgetting the $d$-th point:

$$\pi = \pi_1 \circ \pi_2 \circ \ldots \circ \pi_{2g-1}\,.$$

The formulas for computing $\pi_{d,*}$ of a monomial were already stated in [HM], top of p. 55:

**Formularium.**

a. *Every monomial in the classes $K_i$ $(1 \leq i \leq d)$ and $D_{ij}$ $(1 \leq i < j \leq d)$ on $\mathcal{C}_g^d$ can be rewritten as a monomial $M$ pulled back from $\mathcal{C}_g^{d-1}$ times either a single diagonal $D_{id}$ or a power $K_d^k$ of $K_d$ by a repeated application of the following substitution rules:*

$$\begin{cases} D_{id}D_{jd} \to D_{ij}D_{id} & (i < j < d); \\ D_{id}^2 \to -K_iD_{id} & (i < d); \\ K_dD_{id} \to K_iD_{id} & (i < d). \end{cases}$$

b. *For $M$ a monomial pulled back from $\mathcal{C}_g^{d-1}$:*

$$\begin{cases} \pi_{d,*}(M \cdot D_{id}) = M; \\ \pi_{d,*}(M \cdot K_d^k) = M \cdot \pi^*(\kappa_{k-1}). \end{cases}$$

*Here $\pi : \mathcal{C}_g^{d-1} \to \mathcal{M}_g$ is the forgetful map.*

Note that $\pi_{d,*}(M) = 0$ as it should be, since $\kappa_{-1} = 0$.

Let us discuss some examples, starting in genus 2. On $\mathcal{C}_2^3$ we have:

$$\begin{aligned} 0 &= c_2(\mathbb{F}_3 - \mathbb{E}) \\ &= c_2\big((1+K_1)(1+K_2-D_{12})(1+K_3-D_{13}-D_{23})(1-\lambda_1+\lambda_2)\big) \\ &= (K_1K_2 + K_1K_3 + K_2K_3 - K_1D_{12} - K_1D_{13} - K_1D_{23} \\ &\quad - K_2D_{13} - K_2D_{23} - K_3D_{12} + D_{12}D_{13} + D_{12}D_{23}) \\ &\quad - \lambda_1(K_1 + K_2 + K_3 - D_{12} - D_{13} - D_{23}) + \lambda_2\,. \end{aligned}$$

Upon intersecting this with $D_{13}$ and applying the substitution rules, this becomes:

$$\begin{aligned} 0 &= D_{13} \cdot c_2(\mathbb{F}_3 - \mathbb{E}) \\ &= (2K_1^2 D_{13} + 3K_1K_2D_{13} - 6K_1D_{12}D_{13}) \\ &\quad - \lambda_1(3K_1D_{13} + K_2D_{13} - 2D_{12}D_{13}) + \lambda_2 D_{13}\,. \end{aligned}$$

Pushing this down gives:

$$\begin{aligned} 0 &= \pi_*\big(D_{13} \cdot c_2(\mathbb{F}_3 - \mathbb{E})\big) \\ &= (\pi_1 \circ \pi_2)_*\big(2K_1^2 + 3K_1K_2 - 6K_1D_{12} - \lambda_1(3K_1 + K_2 - 2D_{12}) + \lambda_2\big) \\ &= \pi_{1,*}(3K_1\kappa_0 - 6K_1 - \lambda_1\kappa_0 + 2\lambda_1) \\ &= 3\kappa_0^2 - 6\kappa_0 = 0. \end{aligned}$$

This is no surprise. Instead, pushing-down after intersecting with $D_{12}D_{13}$ gives:

$$\begin{aligned} 0 &= \pi_*\big(D_{12}D_{13} \cdot c_2(\mathbb{F}_3 - \mathbb{E})\big) \\ &= (\pi_1 \circ \pi_2)_*(11K_1^2D_{12} - 6\lambda_1K_1D_{12} + \lambda_2D_{12}) \\ &= \pi_{1,*}(11K_1^2 - 6\lambda_1K_1 + \lambda_2) \\ &= 11\kappa_1 - 6\kappa_0\lambda_1 = 11\kappa_1 - 12\lambda_1 = 10\kappa_1\,, \end{aligned}$$

which implies $\kappa_1 = 0$. Although we knew this already, it nevertheless shows that non-trivial relations can be obtained in this manner. In fact, using $\kappa_1 = 0$, hence $\lambda_1 = 0 = \lambda_2$, we obtain also the relations:

$$\begin{cases} K_1^2 = 0\,; \\ K_1K_2 = 2K_1D_{12}\,; \\ 0 = K_1D_{12} + K_1D_{13} + K_2D_{23} - K_1D_{23} - K_2D_{13} - K_3D_{12} + 2D_{12}D_{13}\,. \end{cases}$$

(The symmetry of the situation was used to obtain the latter relation.) As one can check easily, it follows that $R^*(\mathcal{C}_2^n)$ is Gorenstein (with socle in degree $n$) for $n \leq 3$.

Next we look at genus 3. We know that $\kappa_1 \neq 0$, so there are no relations in $R^1(\mathcal{M}_3)$. The vanishing in codimension 2 resulting from Looijenga's theorem can be made explicit by means of the relations introduced above. First, upon intersecting $c_3(\mathbb{F}_5 - \mathbb{E})$ with the product $D_{12}D_{13}D_{14}D_{15}$ of diagonals and pushing this down to $\mathcal{C}_3$ we find:

$$\begin{aligned} 0 &= c_3\big((1+K)(1+2K)\cdots(1+5K) - \mathbb{E}\big) \\ &= 225K^3 - 85K^2\lambda_1 + 15K\lambda_2 - \lambda_3\,, \end{aligned}$$

which after pushing-down to $\mathcal{M}_3$ gives the relation

$$225\kappa_2 - \frac{55}{8}\kappa_1^2 = 0.$$

Second, intersecting $c_4(\mathbb{F}_5 - \mathbb{E})$ with the product $D_{12}D_{13}D_{45}$ and pushing this down to $\mathcal{M}_3$ one finds (after a calculation) the relation

$$\frac{87}{4}\kappa_1^2 - 162\kappa_2 = 0.$$

Hence we find $\kappa_1^2 = 0 = \kappa_2$, so that $R^*(\mathcal{M}_3) = \mathbb{Q}[\kappa_1]/(\kappa_1^2)$, as was already known to Mumford (cf. [Mu], p. 309).

In genus 4, the interesting codimension is 2: we know that $R^2(\mathcal{M}_4)$ is 1-dimensional; the question is what the precise relation between $\kappa_1^2$ and $\kappa_2$ is. It was determined in [Fa 2]; let us rederive it here. Intersecting $c_4(\mathbb{F}_7 - \mathbb{E})$ with the product $D_{12}D_{13}D_{14}D_{15}D_{67}$ and pushing this down to $\mathcal{M}_4$ we find eventually the relation

$$150\kappa_1^2 - 1600\kappa_2 = 0.$$

Similarly, via $D_{12}D_{13}D_{14}D_{56}D_{57} \cdot c_4(\mathbb{F}_7 - \mathbb{E})$ respectively $D_{12}D_{34}D_{56}D_{57} \cdot c_5(\mathbb{F}_7 - \mathbb{E})$ we find the relations

$$360\kappa_1^2 - 3840\kappa_2 = 0 \qquad \text{and} \qquad -180\kappa_1^2 + 1920\kappa_2 = 0.$$

Lengthy calculations are required to obtain these relations, especially in the last case. It is then reassuring to find that all three are equivalent to

$$\kappa_1^2 = \frac{32}{3}\kappa_2\,,$$

as obtained in [Fa 2] and in accordance with the predicted relation between $\kappa_1^{g-2}$ and $\kappa_{g-2}$ in genus $g$ mentioned before.

This may be a good moment to point out that the relations described in Conjecture 2, geometrically transparent as they are, appear to be rather complicated from a combinatorial point of view. To illustrate this, we observe that $c_g(\mathbb{F}_{2g-1} - \mathbb{E})$ is a polynomial of degree $g$ in roughly $2g^2$ variables. If expanding it completely before applying the substitution and push-down rules of the Formularium were the only option, we would be stuck in genus 5 or 6 already. As we will explain later, there are various ways to overcome this difficulty. For the moment we continue our description of the tautological rings in low genus. All results were obtained by explicitly calculating the relations described in Conjecture 2. We wrote several Maple[1] procedures for the occasion.

Genus 5 gives the first example of a relation not resulting from Looijenga's theorem. It is the relation

$$\kappa_1^2 = \frac{72}{5}\kappa_2\,,$$

establishing that $R^2(\mathcal{M}_5)$ is 1-dimensional, thus providing the first evidence beyond Theorems 1 and 2 for the conjectured Gorenstein property of $R^*(\mathcal{M}_g)$. We also find the predicted relation $\kappa_1^3 = 288\kappa_3$. Let us also give the resulting relations between kappa's and lambda's:

$$\begin{cases} \kappa_1 = 12\lambda_1 \\ \kappa_2 = 10\lambda_1^2 = 20\lambda_2 \\ \kappa_3 = 6\lambda_1^3 = 40\lambda_3 \end{cases}.$$

[1] Maple© is a trademark of the University of Waterloo and Waterloo Maple Software.

Combining these results with the computation of $\kappa_{g-2}\lambda_{g-1}\lambda_g$, we computed the tautological class of the locus of Jacobians of curves of genus 5 in the moduli space $\mathcal{A}_5$ of principally polarized abelian 5-folds. The method is the same as in [Fa 3], §5; see also the concluding remarks. The result is:

$$[\mathcal{J}_5]_Q = \frac{1}{2}[\mathcal{J}_5] = 36\lambda_1^3 - 48\lambda_3 + X$$

in the rational cohomology of a toroidal compactification $\widetilde{\mathcal{A}}_5$, where $X$ satisfies $\lambda_1^{12}X = \lambda_1^9\lambda_3X = 0$.

Note that the relations above, together with Theorem 2, determine the Chow ring of $\mathcal{M}_5$, since Izadi proved that the Chow ring equals the tautological ring in genus 5 ([Iz]).

We remark that for $g \leq 5$

$$R^*(\mathcal{M}_g) = \mathbb{Q}[\kappa_1]/(\kappa_1^{g-1})$$

(and the tautological ring equals the Chow ring). Such a simple description is not available in genus 6 and higher: $\kappa_1^2$ and $\kappa_2$ are independent, as follows from Edidin's result ([Ed]).

We find in genus 6:

$$R^*(\mathcal{M}_6) = \mathbb{Q}[\kappa_1, \kappa_2]/(127\kappa_1^3 - 2304\,\kappa_1\kappa_2\,, 113\,\kappa_1^4 - 36864\,\kappa_2^2).$$

Hence this is still a complete intersection ring. We also have the relations

$$\kappa_3 = \frac{5}{2304}\,\kappa_1^3 \qquad \text{and} \qquad \kappa_4 = \frac{5}{73728}\,\kappa_1^4\,.$$

In genera 7 and 8, the tautological rings are also complete intersection rings, but apparently this is not the case for any genus greater than or equal to 9. As a final example, we give here the tautological ring in genus 9, because it is the first one which is not a complete intersection, and also to give the reader an idea of how complicated these rings become very quickly. The ring $R^*(\mathcal{M}_9)$ is the quotient of $\mathbb{Q}[\kappa_1, \kappa_2, \kappa_3]$ by the ideal generated by

$$\begin{cases} 5195\kappa_1^4 + 3644694\kappa_1\kappa_3 + 749412\kappa_2^2 - 265788\kappa_1^2\kappa_2 \\ 33859814400\kappa_2\kappa_3 - 95311440\kappa_1^3\kappa_2 + 2288539\kappa_1^5 \\ 19151377\kappa_1^5 + 16929907200\kappa_1\kappa_2^2 - 1142345520\kappa_1^3\kappa_2 \\ 1422489600\kappa_3^2 - 983\kappa_1^6 \\ 1185408000\kappa_2^3 - 47543\kappa_1^6 \\ 42019\kappa_1^6 - 1234800\kappa_1^4\kappa_2 \end{cases}$$

(the last generator is actually superfluous). We also have the relations

$$\begin{cases} 1399562496\kappa_4 = 2453760\kappa_2^2 - 65425\kappa_1^4 + 2470320\kappa_1^2\kappa_2 \\ 7223427072\kappa_5 = 307440\kappa_1^3\kappa_2 - 8729\kappa_1^5 \\ 309657600\kappa_6 = \kappa_1^6 \\ 26011238400\kappa_7 = \kappa_1^7 \end{cases}$$

So the dimensions of the vector spaces $R^i(\mathcal{M}_9)$ are $1, 1, 2, 3, 3, 2, 1, 1$ respectively. By straightforward calculations one verifies, first, that the pairings $R^i(\mathcal{M}_9) \times R^{7-i}(\mathcal{M}_9) \to R^7(\mathcal{M}_9) = \mathbb{Q}$ are perfect, i.e., $R^*(\mathcal{M}_9)$ is Gorenstein, second, that $R^*(\mathcal{M}_9)$ satisfies the Hard Lefschetz and Hodge Positivity properties, third, that the proportionalities in degree $g-2=7$ are as conjectured. All the relations above are in the ideal $I_9$ mentioned in Conjecture 2. This completes the proof of Conjectures 1 and 2 in the case $g=9$ at hand.

To prove Conjectures 1 and 2 for other values of $g$, one proceeds entirely analogously. I would like to point out that settling these 2 conjectures can be viewed as a combinatorial problem. For a given value of $g$, the verification of the conjectures amounts to a finite calculation; in principle it can be done on a computer. In practice, we carried this out for $g \leq 15$. As remarked already, this requires an implementation in which the expressions $c_g(\mathbb{F}_{2g-1} - \mathbb{E})$ are *not* expanded completely before the rules of the formularium are applied to them. First, note that every diagonal $D_{ij}$ in $M$ can be used to reduce the number of points by 1, by applying the substitution and push-down rules involving $D_{ij}$ directly to $\mathbb{F}_d$. One finds $d$-pointed jet bundles on $\mathcal{C}_g^k$ (with $k < d$), where the $d$ points appear with multiplicities at the general point of $\mathcal{C}_g^k$. E.g., a $D_{d-1,d}$ occurring in $M$ has the effect of replacing $(1 + K_d - \Delta_d)$ in $c(\mathbb{F}_d)$ by $(1 + 2K_{d-1} - \Delta_{d-1})$; this gives the total Chern class of the $d$-pointed jet bundle on $\mathcal{C}_g^{d-1}$ corresponding to $d$-tuples that contain the $(d-1)$-st point with multiplicity 2. For another example, see the calculation done for $c_3(\mathbb{F}_5 - \mathbb{E})$ in genus 3 above.

Next, the trivial identity

$$c_k(\mathbb{F}_d) = c_k(\mathbb{F}_{d-1}) + (K_d - \Delta_d)\, c_{k-1}(\mathbb{F}_{d-1})$$

can be used to expand the total Chern classes of $d$-pointed jet bundles (possibly with multiplicities) *step by step*. This makes the computations somewhat more manageable. We remark that $\pi_*\big(M \cdot c_j(\mathbb{F}_{2g-1} - \mathbb{E})\big)$ is trivially 0 if not all indices from $\{1, \ldots, 2g-1\}$ are 'covered' by the monomial $M$.

Moreover, we derived a formula expressing $\pi_*\big(\Pi_\alpha \cdot c_k(\mathbb{F}_d)\big)$ directly as a polynomial in the kappa's, where $\alpha = (\alpha_1, \ldots, \alpha_a)$ is a partition of $d$ and $\Pi_\alpha$ is a product of $d-a$ diagonals corresponding to $\alpha$:

$$\Pi_\alpha = D_{1,2}D_{1,3}\cdots D_{1,\alpha_1}D_{\alpha_1+1,\alpha_1+2}\cdots D_{\alpha_1+1,\alpha_1+\alpha_2}\cdots D_{d-\alpha_a+1,d-\alpha_a+2}\cdots D_{d-\alpha_a+1,d} \quad .$$

In joint work with Zagier, this formula was rewritten in a form involving formal power series. This form was then used to prove the following statements about the relations in the ideal $I_g$ introduced in Conjecture 2:

a. $I_g$ contains no relations in codimensions $\leq g/3$.

b. For $g$ of the form $3k-1$ with $k$ an integer, there is a unique relation in codimension $k$ in $I_g$.

c. Write $g = 3k - 1 - \ell$, with $k$ and $\ell$ positive integers. There exists an upper bound for the number of relations in $I_g$ in codimension $k$, which only depends on $\ell$.

The first statement is consistent with the second half of part (b) of Conjecture 1 (but doesn't prove it). Actually, the improved stability result of Harer ([Ha 2]) essentially implies that there are no relations between the kappa's in codimensions $\leq g/3$. (I learned of this result after formulating Conjecture 1.)

The unique relation in codimension $k$ and genus $3k-1$ is given as follows. Define rational numbers $a_i$ for $i \geq 1$ via the identity

$$\exp\left(-\sum_{i=1}^{\infty} a_i t^i\right) = \sum_{n=0}^{\infty} \frac{(6n)!}{(2n)!(3n)!} t^n$$

of formal power series. Then the said relation is the coefficient of $t^k$ in

$$\exp\left(\sum_{i=1}^{\infty} a_i \kappa_i t^i\right) .$$

As to the third statement, the available computational evidence suggests that the *actual* number of relations in $R^*(\mathcal{M}_g)$ in codimension $k$ and genus $g = 3k-1-\ell$ depends only on $\ell$, whenever $2k \leq g-2$ (i.e., $k \geq \ell + 3$). Assuming this, and denoting this number by $a(\ell)$, we know that $a(\ell) = 1, 1, 2, 3$ for $\ell = 0, 1, 2, 3$ respectively; further computations assuming Conjecture 1 give the following results:

| $\ell$ | 0 | 1 | 2 | 3 | 4 | 5 | 6 | 7 | 8 | 9 |
|---|---|---|---|---|---|---|---|---|---|---|
| $a(\ell)$ | 1 | 1 | 2 | 3 | 5 | 6 | 10 | 13 | 18 | 24 |

Zagier and I have a favourite guess as to what the function $a$ might be (but there are many functions with 10 prescribed values).

Detailed results concerning the formula expressing $\pi_*\big(\Pi_\alpha \cdot c_k(\mathbb{F}_d)\big)$ as a polynomial in the kappa's as well as the work with Zagier will appear in [Fa 4] and [FZ].

So far, we have only considered relations resulting from the triviality that divisors of negative degree on a curve are not effective (or rather from the dual statement). Naturally, many more relations can be produced with the same method. Start with a triple $(g, d, r)$ such that the locus $\{\operatorname{rank} \varphi_d \leq d - r\}$ has the expected codimension, so that its class can be computed using Porteous's formula. If the fibres of the map to $\mathcal{M}_g$ have the expected dimension $r$, then cutting the locus with $r$ sufficiently general divisors will give a locus that maps finitely onto the locus in $\mathcal{M}_g$ of curves possessing a $g^r_d$. Often, divisors can be chosen whose classes lie in the tautological ring of $\mathcal{C}^d_g$, and often there are quite a few possible choices for such divisors. Every such choice leads to a formula for the class of the locus in $\mathcal{M}_g$ of curves possessing a $g^r_d$, with a certain multiplicity; if the multiplicity can be computed—this is often the case—we obtain the actual class, as an element of the tautological ring. (These observations are certainly not new; see for instance the papers [Di 1], [HM], [Mu], [Ra].)

Equating the various formulas for this class leads then to relations in the tautological ring. However, if one is after such relations (we are), it is a lot easier to cut with *fewer* than $r$ divisors whose classes lie in the tautological ring; pushing-down the resulting class to $\mathcal{M}_g$ gives a polynomial in the kappa's which is 0 in the ring, since the fibres of the map to $\mathcal{M}_g$ are positive-dimensional. In this way we obtain relations without having to compute multiplicities.

Examples:

a. 1-dimensional linear systems. The loci involved have the expected dimensions. The locus of divisors moving in a $g^1_d$ has dimension $2g-4+2d$. Cutting its class with $K_1$ (i.e., requiring that the first point be in a fixed canonical divisor) gives a class that pushes down to the class of curves with a $g^1_d$, with multiplicity $(d-1)!(2g-2)$. Cutting it with $D_{12}$ instead, we find that class with multiplicity $(d-2)!(2g+2d-2)$. The 2 formulas are worked out in [Ra], §5, for $d=2$ resp. 3. In [Fa 2] we proved that the resulting relation expresses $\kappa_{g-2}$ (resp. $\kappa_{g-4}$) in lower kappa's, for $g \geq 4$ (resp. for $g \geq 7$). The other option is to push-down directly. This produces a non-trivial relation in codimension $g-2d+1$. It appears that together these relations are sufficient to prove the first half of part (b) of Conjecture 1, but I have not been able to carry out the calculation.

b. Plane quintics. They form a 12-dimensional subvariety of $\mathcal{M}_6$. Inside $\mathcal{C}^5_6$ the divisors moving in a $g^2_5$ form a 14-dimensional locus. Note that it has 2 components: one lies over the plane quintics; the other lies over the hyperelliptic locus, with 3-dimensional fibres (the $g^2_5$ equals the $g^2_4$ plus an arbitrary base point). Cutting the class with 2 divisors and pushing-down, we find the class of the plane quintics in $\mathcal{M}_6$, with a certain multiplicity. There are 5 ways of choosing 2 tautological divisors: $D_{12}D_{34}$, $D_{12}D_{13}$, $K_1D_{23}$, $K_1D_{12}$ and $K_1K_2$. As one checks easily, these give the class of the plane quintics with multiplicities 240, 90, 360, 60 and 600 respectively. (E.g., $K_1D_{23}$ puts the first point in a fixed canonical divisor; the second and third point coincide, so this is a point of tangency on a line through the first point; this fixes the fourth and fifth point, but not their order, for a total of $10 \cdot 18 \cdot 2 = 360$ possibilities.) The 5 formulas for the class give 4 relations between $\kappa_1^3$, $\kappa_1\kappa_2$ and $\kappa_3$. It turns out these relations have rank 2 (the actual rank, as we know now). (When I did these calculations originally (Summer 1991), the result strongly suggested to me that $R^3(\mathcal{M}_6)$ should be 1-dimensional; this was instrumental in formulating the first version of Conjecture 1 later that year.) The class of the plane quintics turns out to be $\frac{35}{3072}\kappa_1^3$.

c. Hyperelliptic curves. If Porteous's formula can be applied to compute the class of divisors moving in a $g^r_d$, it can also be applied to compute the class of the divisors moving in the dual system (a $g^{g-d+r-1}_{2g-2-d}$). Let us apply this to hyperelliptic curves. The class of divisors moving in a $g^{g-2}_{2g-4}$ is

$$\Delta_{g-2,2}\big(c(\mathbb{F}_{2g-4}-\mathbb{E})\big) = (c^2_{g-2} - c_{g-1}c_{g-3})(\mathbb{F}_{2g-4}-\mathbb{E}).$$

This has $(g-2)$-dimensional fibres over the hyperelliptic locus. Cutting it with $c < g-2$ divisors and pushing-down, we find relations in codimension $c$. There are many possible choices for the $c$ divisors; e.g., one can take $c$ diagonals, in which case the possibilities correspond to the partitions of $c$. It appears that these relations (together with the relation mentioned in (a) expressing $\kappa_{g-2}$ in lower kappa's) generate the entire ideal of relations between the kappa's, but the computations are quite cumbersome.

Concluding remarks:

a. As we saw before, Looijenga's theorem (Theorem 1) in degree $g-1$ implies Diaz's upper bound $g-2$ for the dimension of a complete subvariety of $\mathcal{M}_g$. One may view Theorem 2 as indicating that there is no intersection-theoretical obstruction to the existence of a $(g-2)$-dimensional complete subvariety. Thinking along these lines, Theorem 1 in degree $g-2$ as well as Conjecture 1 put severe constraints on the possible $(g-2)$-dimensional complete subvarieties, whereas complete subvarieties of dimensions $\leq g/3$ are unconstrained from this point of view. In the known existence results, the genus is exponential in the dimension of the complete subvariety; for all $g \geq 3$, complete curves exist; $\mathcal{M}_8$ contains complete surfaces. Perhaps (for now) $\mathcal{M}_6$ and $\mathcal{M}_7$ are more natural places to look for complete surfaces than $\mathcal{M}_4$ and $\mathcal{M}_5$.

b. In degree $g-2$, we have found 'experimentally' several explicit proportionalities. Even deducing them from part (c) of Conjecture 1 appears to be non-trivial, so we are far from proving them. We state the most relevant ones:

$$[\mathcal{H}_g]_Q = \frac{1}{2}\,[\mathcal{H}_g] = \frac{(2^{2g}-1)\,2^{g-2}}{(2g+1)(g+1)!}\,\kappa_{g-2}\,;$$
$$\lambda_{g-2} = |B_{2g-2}|\,\frac{(2g-1)\,2^{g-1}}{(2g-2)(g-1)!}\,\kappa_{g-2}\,.$$

Here $\mathcal{H}_g$ is the hyperelliptic locus. The formula for its class cries out for a direct proof. The second formula, together with Theorem 2, leads to

$$\lambda_{g-1}^3 = \frac{|B_{2g-2}B_{2g}|}{(2g-2)(2g)}\,\frac{1}{(2g-2)!}\quad.$$

This is the contribution from the constant maps, as it occurs in the theory of counting curves of higher genus on threefolds (see [BCOV], §5.13, (5.54)). Before, this number was known only for $g \leq 4$ (cf. [Fa 3]) and no conjectural formula was known. Needless to say, the 3 formulas above are true for all $g \leq 15$.

*Acknowledgements.* I would like to thank Robbert Dijkgraaf, Bill Fulton, Gerard van der Geer, Eduard Looijenga, Shigeyuki Morita, Ragni Piene, Piotr Pragacz, Michael Thaddeus, Chris Zaal and Don Zagier. This research has been made possible by a fellowship of the Royal Netherlands Academy of Arts and Sciences.

Faculteit WINS, Universiteit van Amsterdam, Plantage Muidergracht 24,
1018 TV Amsterdam, The Netherlands

September 1, 1996 — May 31, 1997:
Institut Mittag-Leffler, Auravägen 17, S–182 62 Djursholm, Sweden

e-mail: faber@fwi.uva.nl